\documentclass[12pt]{amsart}

\usepackage{amsmath}
\usepackage{amscd}
\usepackage{amssymb}

\newtheorem{theorem}{Theorem}[section]
\newtheorem{lemma}[theorem]{Lemma}

\theoremstyle{definition}
\newtheorem{definition}[theorem]{Definition}

\theoremstyle{remark}
\newtheorem{remark}[theorem]{Remark}

\numberwithin{equation}{section}

\input{prepictex}
\input{pictex}
\input{postpictex}
\chardef\bslash=`\\ 

\makeatletter
\def\verbatim{\interlinepenalty\@M \@verbatim
  \leftskip\@totalleftmargin\advance\leftskip2pc
  \frenchspacing\@vobeyspaces \@xverbatim}
\makeatother
\newcounter{picture}

\newcommand{\cl}[1]{{\mathcal{#1}}}
\newcommand{\bb}[1]{{\mathbb{#1}}}
\newcommand{\fk}[1]{{\mathfrak{#1}}}

\newcommand{\id}{{\bf 1}}

\newcommand{\PGL}{{\text{\rm{PGL}}}}

\newcommand{\aut}{{\text{\rm Aut}}}

\begin{document}

\title{Boundary $C^*$-algebras for acylindrical groups}

\author{Guyan Robertson}
\address{School of Mathematics and Statistics, University of Newcastle, NE1 7RU, U.K.}
\email{a.g.robertson@newcastle.ac.uk}

\subjclass[2000]{Primary 20E08, 46L80}

\date{}

\keywords{Acylindrical group, boundary, Cuntz-Krieger algebra}

\begin{abstract} Let $\Delta$ be an infinite, locally finite tree with more than two ends. Let $\Gamma<\aut(\Delta)$ be an acylindrical uniform lattice. Then the boundary algebra $\cl A_\Gamma = C(\partial\Delta)\rtimes \Gamma$ is a simple Cuntz-Krieger algebra whose K-theory is determined explicitly.
\end{abstract}

\maketitle

\section{Introduction}

Let $\Delta$ be an infinite, locally finite tree with more than two ends and with boundary $\partial\Delta$.
Let $k$ be a positive integer and let $\Gamma$ be a group of automorphisms of $\Delta$ without inversion and with no proper invariant subtree.
Say that $\Gamma$ is $k$-acylindrical if the stabilizer of any path of length $k$ in $\Delta$ is trivial \cite{bp, sel}.
The group $\Gamma$ is {\it acylindrical} if it is $k$-acylindrical for some integer $k\ge 1$.
The main result of this article is
\begin{theorem}\label{main0}
Let $\Gamma<\aut(\Delta)$ be an acylindrical uniform lattice. Then the boundary algebra $\cl A_\Gamma = C(\partial\Delta)\rtimes \Gamma$ is a simple Cuntz-Krieger algebra.
\end{theorem}

The action of $\Gamma$ on $\partial\Delta$ is amenable, so the maximal crossed product $\cl A_\Gamma = C(\partial\Delta)\rtimes \Gamma$ coincides with the reduced crossed
product and is nuclear by Proposition 4.8 and Th\'eor\`eme 4.5 in \cite{ad}.
The algebra $\cl A_\Gamma$ is described in Section \ref{cp} below.
Cuntz-Krieger algebras were introduced in \cite{ck} and are classified up to isomorphism by their K-theory \cite{k}. A special case, where $\Gamma$ is a free uniform tree lattice, was studied in \cite{rtree} by different methods.
  The K-groups of the boundary algebra $\cl A_\Gamma$ are isomorphic to the Bowen-Franks invariants of flow equivalence for a certain subshift of finite type associated with the geodesic flow \cite{c}. This subshift was studied in \cite[6.3]{bp}.

The K-groups of the algebra $\cl A_\Gamma$ may be computed explicitly. For example, if  $\Gamma=\bb Z_{l+1}*\bb Z_{m+1}$ acts on its Bass-Serre tree, where $l,m \ge 1$, then
$K_0(\cl A_\Gamma)=\bb Z_{lm-1}$ via an isomorphism sending the class of the identity idempotent to $l+1$.
It follows that $\cl A_\Gamma \cong M_{l+1}\otimes \cl O_{lm}$, where $\cl O_n$ denotes the Cuntz algebra, which is generated by $n$ isometries on a Hilbert space whose range projections sum to the identity operator \cite{c0}.

\begin{remark}\label{qi}
  The algebra $\cl A_\Gamma$ depends only on $\Gamma$. The group $\Gamma$ and the tree $\Delta$ are quasi-isometric; more precisely, for any base vertex in $\Delta$, the natural mapping from $\Gamma$ onto the orbit of this vertex is a quasi-isometry from $\Gamma$ to $\Delta$. This is a special case of the ``Fundamental Observation of Geometric Group Theory'' \cite[Theorem~IV.23]{har}. The mapping from $\Gamma$ to $\Delta$ induces a $\Gamma$-equivariant homeomorphism of the boundaries of $\Gamma$ and $\Delta$.
\end{remark}

\section{Background}

The edges of the tree $\Delta$ are directed, and each geometric edge of $\Delta$ corresponds to two directed edges. Let $\Delta ^0$ denote the set of vertices and $\Delta^1$  the set of directed edges of $\Delta$. There is a distance function $d$ defined on the geometric realization of $\Delta$ which assigns unit length to each edge. Choose an orientation on the set of edges which is invariant under $\Gamma$. This orientation consists of a partition of $\Delta^1$ and a bijective involution
$$e\mapsto \overline e : \Delta^1 \to \Delta^1$$
which interchanges the two components of $\Delta^1$. Each directed edge $e$ has an initial vertex $o(e)$ and a terminal vertex $t(e)$ such that $o(\overline e)=t(e)$.

Let $\Gamma$ be a group of automorphisms of $\Delta$ without inversion and with no proper invariant subtree.
Say that $\Gamma$ is $k$-acylindrical, where $k\ge 1$, if the stabilizer of any path of length $k$ in $\Delta$ is trivial \cite{bp}.
(In \cite{sel} such a group $\Gamma$ is said to be $(k-1)$-acylindrical.)
To say that $\Gamma$ is $1$-acylindrical is the same as saying that $\Gamma$ acts freely on $\Delta^1$. For example, the action of a free product $\Gamma=\Gamma_1*\Gamma_2$ on the associated Bass-Serre tree \cite[I.4.1]{ser} is $1$-acylindrical. If $\Gamma=\Gamma_1*_{\Gamma_0}\Gamma_2$ is a free product with amalgamation over $\Gamma_0$, then the action of $\Gamma$ on the associated Bass-Serre tree  is $2$-acylindrical if $\Gamma_0$ is malnormal in each $\Gamma_j$, i.e. $g^{-1}\Gamma_0g \cap \Gamma_0 = \{ 1 \}$ for all $g \in \Gamma_j - \Gamma_0$, $j=1,2$. Every small splitting of a torsion free hyperbolic group gives rise to a 3-acylindrical action \cite{sel}.

The boundary $\partial\Delta$ is the set of equivalence classes of infinite semi-geodesics in $\Delta$, where two semi-geodesics are said to be equivalent if they agree except on finitely many edges. For the rest of this article we make the following assumptions.

\medskip
\noindent\textbf{Standing Hypotheses}
\begin{description}
\item[(1)] $\Delta$ is an infinite locally finite tree with more than two boundary points.
\item[(2)] $\Gamma<\aut(\Delta)$ is a uniform tree lattice.
\item[(3)] $\Gamma$ acts without inversion and with no proper invariant subtree.
\item[(4)] $\Gamma$ is $k$-acylindrical, where $k\ge 1$.
\end{description}

\begin{remark}
The standing hypotheses imply that $\partial\Delta$ is uncountable.
\end{remark}
\begin{remark}
  The assumption that $\Gamma$ has no proper invariant subtree is part of the definition of ``acylindrical'' in \cite{sel} (but not in \cite{bp}, which is why it is emphasised separately here). It implies that the action of $\Gamma$ on $\partial\Delta$ is minimal. The assumption could have been omitted here if $\partial\Delta$ were replaced throughout by the limit set $\Lambda_\Gamma$.
\end{remark}
\begin{remark}
  The assumption that $\Gamma$ is a \textit{uniform} lattice is natural and necessary in the context of the theory of \cite{rs1}. It would be interesting to study $\cl A_\Gamma$ for non-uniform tree lattices. For example, \cite{bp} provides an example of a 5-acylindrical action of the Nagao lattice $\PGL_2(\bb F_q[t])$.
\end{remark}

\section{Cuntz-Krieger algebras}\label{HRCK}

It is convenient to use the approach to Cuntz-Krieger algebras developed in \cite{rs1}.
Choose a nonzero matrix  $M$ with entries in $\{0,1\}$. For $m\ge 0$, let $W_m$ denote the set of all words of length $m+1$ based on the alphabet $A$ and the transition matrix $M$. A word $w\in W_m$ is a formal product $w=a_0a_1\dots a_m$, where $a_j\in A$ and
$M(a_{j+1}, a_j)=1$, $0\le j\le m-1$. Define $o(w) = a_0$ and $t(w)=a_m$.

Fix a nonempty finite or countable set $D$ (whose elements are ``decorations'')
and a map $\delta : D \to A$.
Let $\overline W_m = \{ (d,w) \in D \times W_m ;\ o(w) = \delta (d) \}$, the set of
``decorated words'' of length $m+1$, and identify
$D$ with $\overline W_0$ via the map $d \mapsto (d,\delta(d))$.
Let $W = \bigcup_m W_m$ and  $\overline W = \bigcup_m \overline W_m$, the sets of all
words and all decorated words respectively.
Define $o: \overline W_m \to D$ and $t: \overline W_m \to A$ by
$o(d,w) = d$ and $t(d,w) = t(w)$.

Let $u=a_0a_1\dots a_m\in W_m$ and $v=b_0b_1\dots b_n \in W_n$. If $t(u) =o(v)$, then there exists a unique
product $uv\in W_{m+n}$ defined by
$$uv=a_0a_1\dots a_mb_1\dots b_n.$$

If $w=a_0a_1\dots a_l\in W_l$ where $l \ge 0$ and if $p \ne 0$, we say that
$w$ is {\em $p$-periodic} if $a_{j+p}=a_j$ whenever both sides are defined.
Assume that the nonzero $\{0,1\}$-matrix $M$ has been chosen so that the following conditions from \cite{rs1} hold.

\begin{description}
\item[(H2)] If $a, b \in A$ then there exists $w\in W$ such that $o(w)=a$ and $t(w)=b$.
\item[(H3)] For each nonzero integer $p$, there exists some  $w \in W$ which is not $p$-periodic.
\end{description}

\begin{definition}\label{CKdefinition} \cite{rs1}
The $C^*$-algebra $\cl A_D=\cl A (A,D,M)$ is the universal $C^*$-algebra
generated by a family of partial isometries
$\{s_{u,v};\ u,v \in \overline W \ \text{and} \ t(u) = t(v) \}$
satisfying the relations
\begin{subequations}\label{rel1*}
\begin{eqnarray}
{s_{u,v}}^* &=& s_{v,u} \label{rel1a*}\\
s_{u,v}s_{v,w}&=&s_{u,w} \label{rel1b*}\\
s_{u,v}&=&\displaystyle\sum_
{\substack{w\in W_1,\\
o(w)=t(u)=t(v)}}
s_{uw,vw}
\label{rel1c*}\\
s_{u,u}s_{v,v}&=&0 ,\ \text{for} \ u,v \in \overline W_0, u \ne v. \label{rel1d*}
\end{eqnarray}
\end{subequations}
\end{definition}

\begin{remark}\label{id}
If $D$ is finite (as is the case in this article), then $\cl A_D$ is isomorphic to a simple Cuntz-Krieger algebra with identity element $\id = \displaystyle\sum_{u\in \overline W_0}s_{u,u}$ \cite{rws}.
\end{remark}

\begin{remark}\label{decorate}
Two decorations $\delta_1 : D_1 \to A$ and $\delta_2 : D_2 \to A$ are said to be {\it equivalent} \cite[Section 5]{rs1} if there is a bijection $\eta: D_1 \to D_2$ such that $\delta_1 = \delta_2 \eta$. Equivalent decorations $\delta_1, \delta_2$ give rise to isomorphic algebras $\cl A_{D_1}, \cl A_{D_2}$.
\end{remark}

\begin{remark}\label{OA}
  Denote by $\cl A_A$ the Cuntz-Krieger algebra with decorating set $A$ and with $\delta$ the identity map. The algebra $\cl A_A$ is isomorphic to the algebra $\cl O_{M^t}$ generated by
a set of partial isometries $\{S_a ; a\in A\}$ satisfying the relations
$S_a^*S_a=\sum_bM(b,a)S_bS_b^*$ \cite[Remark 3.11]{rs1}.
If $A$ contains $n$ elements and $M(b, a)=1$, for all $a,b\in A$, then $\cl A_A$ is the Cuntz algebra $\cl O_n$ generated by $n$ isometries whose range projections sum to the identity operator \cite{c0}.
\end{remark}

\section{The algebra associated with an acylindrical group}\label{cp}

The Standing Hypotheses (1)--(4) are now in force. A \emph{geodesic} $\gamma$ in $\Delta$  is a sequence
$(s_j)_{j=-\infty}^\infty$ of vertices such that $d(s_i,s_j)=|i-j|$. A \emph{directed segment} $\sigma$ of length $n$ is a sequence $(s_0, s_1, \dots, s_n)$ of vertices such that $d(s_i, s_j) = |i-j|$. Denote such a directed segment by $[s_0, s_n]$ and let $\fk S_n$ be the set of directed segments of length $n$ in $\Delta$. Since the group $\Gamma$ is $k$-acylindrical, $\Gamma$ acts freely on the set $\fk S_k$.

The alphabet $A$ is defined to be $\Gamma\backslash\fk S_{k+1}$, the set of $\Gamma$-orbits of directed segments of length $k+1$ in $\Delta$.
Since $\Gamma$ is a uniform lattice, $A$ is finite.

\begin{figure}[htbp]
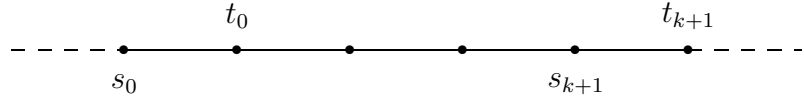
\label{transition}
\hfil
\centerline{
\beginpicture
\setcoordinatesystem units <1.5cm, 1.5cm>
\setplotarea  x from -4 to 4,  y from -0.3 to 0.3
\multiput {$_\bullet$} at -2 0 *5 1 0 /
\putrule from -2 0 to 3 0
\setdashes
\putrule from 3 0 to 4 0
\putrule from -3 0 to -2 0
\put {$s_0$}    at -2 -.3
\put {$s_{k+1}$}    at 2 -.3
\put {$t_0$}    at -1 .3
\put {$t_{k+1}$}    at 3 .3
\endpicture
}
\hfil
\caption{The condition $M(b,a)=1$.}
\end{figure}

Define a matrix $M$ with entries in $\{0,1\}$ as follows.
If $a,b \in A$, we say that $M(b,a)=1$ if and only if $a=\Gamma \sigma$ and $b=\Gamma \tau$, where  $\sigma=(s_0, s_1, \dots, s_{k+1})$, $\tau=(t_0, t_1, \dots, t_{k+1})$ are directed segments such that $s_{j+1}=t_j$, $0\le j \le k$.
The definition is illustrated in Figure~\ref{transition}.

\refstepcounter{picture}

As in Section \ref{HRCK}, $W_m$ denotes the set of all words of length $m+1$ based on the alphabet $A$ and transition matrix $M$.
Let $\fk W_m = \Gamma\backslash \fk S_{m+k+1}$ and let $\fk W = \bigcup_m \fk W_m$. There is a map
$$\alpha: \fk W_m \to W_m$$
defined by
$$
\alpha(\Gamma (s_0, s_1, \dots, s_{m+k+1})) = (\Gamma [s_0, s_{k+1}])(\Gamma [s_1, s_{k+2}])\dots (\Gamma [s_m, s_{m+k+1}]).
$$

\begin{lemma}\label{L1}
The map $\alpha$ is a bijection from $\fk W_m$ onto $W_m$.
\end{lemma}

\begin{proof}
Suppose that $\alpha(\Gamma \sigma)=\alpha(\Gamma \tau)$, where $\sigma=(s_0, s_1, \dots, s_{m+k+1})$ and  $\tau=(t_0, t_1, \dots, t_{m+k+1})$. Then $\Gamma [s_j, s_{j+k+1}] = \Gamma [t_j, t_{j+k+1}]$, $0\le j \le m$.
For each $j$ there exists $g_j\in \Gamma$ such that $g_j [s_j, s_{j+k+1}] = [t_j, t_{j+k+1}]$, and $g_j$ is uniquely determined, since $\Gamma$ acts freely on the set of segments of length $k$.
Now if $0\le j \le m-1$, then $[s_j, s_{j+k+1}]\cap [s_{j+1}, s_{j+k+2}]= [s_{j+1}, s_{j+k+1}]$ is a segment of length $k$.
Also $g_j [s_{j+1}, s_{j+k+1}] = [t_{j+1}, t_{j+k+1}]=g_{j+1} [s_{j+1}, s_{j+k+1}]$. Therefore
$g_j=g_{j+1}$. It follows that $g_0\sigma=\tau$ and so $\Gamma \sigma=\Gamma \tau$. This proves injectivity.

To prove surjectivity, suppose that $w=a_0a_1\dots a_m\in W_m$. Then by definition $a_j= \Gamma \sigma_j$ for $0\le j \le m$, where $\sigma_j\cap\sigma_{j+1}$ is a segment of length $k$, for $1\le j \le m-1$.
It follows that there is a directed segment $\sigma =(s_0, s_1, \dots , s_{m+k+1})\in \fk S_{m+k+1}$ such that
$\sigma_j=[s_j, s_{j+k+1}], 0\le j \le m$. Thus $\alpha (\Gamma\sigma)=w$.
\end{proof}

Now fix a vertex $P\in\Delta^0$. Let $\overline{\fk W}_m$ denote the set of directed segments of length $m+k+1$ which begin at $P$ and let $\overline{\fk W} = \bigcup_{m \ge 0} \overline{\fk W}_m$. The decorating set is $D=\overline{\fk W}_0$, and the decorating map $\delta : D \to A$ is defined by $\delta(d) = \Gamma d$.
Define $\overline{\alpha} : \overline{\fk W} \to \overline{W}$ by
$$\overline{\alpha}(\sigma)=\left(o(\sigma), \alpha(\Gamma \sigma)\right),$$
where
$o(\sigma)=(s_0, s_1, \dots , s_{k+1})$
 is the initial segment of length $k+1$ of $\sigma$. Also define $t(\sigma)$ to be the
final segment of length $k+1$ of $\sigma$.

\begin{lemma} \label{L2} The map $\overline{\alpha}$ is a bijection from $\overline{\fk W}_m$ onto $\overline{W}_m$, for each $m \ge 0$.
\end{lemma}

\begin{proof}
If $\overline{\alpha}(\sigma_1)=\overline{\alpha}(\sigma_2)$, then $o(\sigma_1)=o(\sigma_2)$; moreover $\Gamma \sigma_1=\Gamma \sigma_2$ by Lemma \ref{L1}.  Since $\Gamma$ acts freely on $\fk S_k$, it
follows that $\sigma_1=\sigma_2$.  Therefore $\overline{\alpha}$ is injective.

To see that $\overline{\alpha}$ is surjective, let $\overline{w} = (d,w) \in \overline{W}_m$, where $w
\in W_m$ and $d \in D$. By Lemma \ref{L1}, there exists $\sigma \in \fk S_{m+k+1}$ such that $\alpha( \Gamma \sigma)=w$. Now
\begin{equation*}
\Gamma d=\delta (d)=o(w)=o(\alpha(\Gamma \sigma))= \Gamma o(\sigma).
\end{equation*}
Replacing $\sigma$ by $g \sigma$ for suitable $g \in \Gamma$ ensures that $o(\sigma)=d$.
Then $\sigma\in \overline{\fk W}_m$ and $\overline{\alpha} (\sigma)=\overline{w}$.
\end{proof}

Recall that each $\omega\in\partial\Delta$ is represented by a unique semi-geodesic $[s_0,\omega)$ with initial vertex $s_0\in \Delta^0$.
If $\sigma$ is a directed segment with initial vertex $t_0$, let
$$
\Omega(\sigma) = \left \{ \omega \in \partial\Delta : [t_0,\omega) \ \text{ contains }  \sigma \right \}.
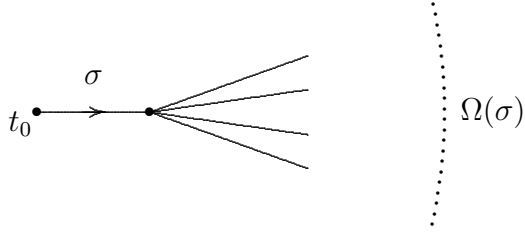
$$
\centerline{
\beginpicture
\setcoordinatesystem units <1.5cm, 1.5cm>    
\setplotarea  x from 0 to 4,  y from -1.5 to 1.5
\put{${\sigma}$}   [b] at 0.5 0.25
\put{$_\bullet$}  at  0  0
\put{$t_0$}  at  -0.15 -0.1
\put{$_\bullet$}  at  1 0
\put{${\Omega(\sigma)}$}   [l] at 3.75 0
\arrow <6pt> [.2, .67] from  0.4 0  to 0.6 0
\setlinear
\plot    0  0  1 0  2.4 0.5 /
\plot    1  0  2.4 0.2 /
\plot    1  0  2.4 -0.2 /
\plot    1  0  2.4 -0.5 /
\setplotsymbol({$\cdot$}) \plotsymbolspacing=4pt
\circulararc 25 degrees from 3.5 -1 center at -1 0
\endpicture
}
The boundary $\partial\Delta$ has a natural compact totally disconnected topology generated by sets of the form $\Omega(\sigma)$ where $\sigma\in \overline{\fk W}$ \cite[I.2.2]{ser}.
The group $\Gamma$ acts on $\partial\Delta$, and one can form the crossed product $C^*$-algebra $C(\partial\Delta) \rtimes \Gamma$.
This is the universal $C^*$-algebra  generated by the commutative $C^*$-algebra $C(\partial\Delta)$ and the image of a unitary
representation $\pi$ of $\Gamma$, satisfying the covariance relation
\begin{equation}\label{cov0}
f(g^{-1}\omega) = \pi(g)\cdot f \cdot \pi(g)^{-1}(\omega)
\end{equation} for
$f \in C(\partial\Delta)$, $g \in \Gamma$ and $\omega\in\partial\Delta$.
It is convenient to denote $\pi(g)$ simply by $g$. Equation (\ref{cov0}) implies that for each clopen set $E\subset\partial\Delta$,
\begin{equation}\label{cov}
\chi_{gE} = g\cdot \chi_E\cdot  g^{-1}.
\end{equation}
The indicator function $\chi_E$ is continuous and is regarded as an element of the crossed product algebra via the embedding $C(\partial\Delta)\subset C(\partial\Delta) \rtimes \Gamma$. The following is a more precise version of Theorem \ref{main0}.

\begin{theorem}\label{main1}
Let $\cl A_\Gamma = C({\partial \Delta}) \rtimes \Gamma$. Then $\cl A_\Gamma$ is isomorphic to the Cuntz-Krieger algebra $\cl A_D$ associated with the alphabet $A$, the decorating set $D$ and the transition matrix $M$.
\end{theorem}

\begin{proof}
The isomorphism $\phi: \cl A_D \to C(\partial \Delta)\rtimes \Gamma$ is defined as follows.
Let $\overline w_j =\overline\alpha(\sigma_j) \in \overline W$, $j=1,2$, with $t(\overline w_1)=t(\overline w_2)$. By the definition of $\overline\alpha$, there is an element $g\in \Gamma$
such that $gt(\sigma_1)=t(\sigma_2)$. Recall that $t(\sigma_1)$ and $t(\sigma_2)$ are segments of length $k+1$. Therefore $g$ is unique, since $\Gamma$ acts freely on $\fk S_k$. Define the homomorphism $\phi$ by
$$
\phi(s_{\overline w_2,\overline w_1})= g\chi_{\Omega(\sigma_1)}=\chi_{\Omega(\sigma_2)}g.
$$
This equation defines a $*$-homomorphism of $\cl A_\Gamma$ because the operators of the form
$\phi(s_{\overline w_2,\overline w_1})$ are easily seen to satisfy the relations (\ref{rel1*}).
Since the algebra $\cl A_D$ is simple \cite[Theorem 5.9]{rs1}, $\phi$ is injective.

Now $\chi_{\Omega(\sigma)}= \phi(s_{\overline w, \overline w})$,
where $\sigma\in \overline{\fk W}$ and $\overline w =\overline\alpha(\sigma)$. Since the sets $\Omega(\sigma)$, $\sigma\in \overline{\fk W}$,  form a basis for the topology of $\Omega$, the linear span of $\{ \chi _{\Omega (\sigma)} ; \sigma\in \overline{\fk W}\}$ is dense in $C(\partial\Delta)$.
It follows that the range of $\phi$ contains $C(\partial\Delta)$.
To show that $\phi$ is surjective, it therefore suffices to show that the range of $\phi$ contains $\Gamma$.

Let $g\in \Gamma$ and choose an integer $m\ge d(P, g^{-1}P)$. Let $\sigma\in \overline{\fk W}_m$. Then $\sigma$ is a
directed segment of length $m+k+1$ with initial vertex $P$ and final vertex $Q$, say. Let $\sigma''$ be the directed segment with initial vertex $g^{-1}P$ and final vertex $Q$. Since $m\ge d(P, g^{-1}P)$, it follows that $t(\sigma'')=t(\sigma)$.

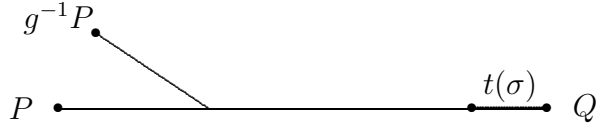
\centerline{
\beginpicture
\setcoordinatesystem units <1cm,1cm>   
\setplotarea x from 0 to 6.5, y from -0.5  to 1.7         
\putrule from 0 0 to 6.5 0
\setlinear
\plot  0.5 1  2 0 /
\put{$P$}  at  -0.5 0
\put{$Q$}  at  7 0
\put{$g^{-1}P$}  at  0 1.2
\put{$_\bullet$}  at  0 0
\put{$_\bullet$}  at  0.5  1
\put{$_\bullet$}  at  5.5 0
\put{$_\bullet$}  at  6.5 0
\put{$t(\sigma)$}  at  6 0.3
\setplotsymbol({$\cdot$}) \plotsymbolspacing=1pt
\plot  5.5 0   6.5 0 /
\endpicture
}
Let $\sigma'=g\sigma''$. Then $\sigma'$ is a path beginning at $P$ and $t(\sigma') = gt(\sigma)$. Let $\overline w_1 = \overline \alpha(\sigma)$ and $\overline w_2 = \overline \alpha(\sigma')$.
Then $t(\overline w_1)=t(\overline w_2)$ and
$g\chi_{\Omega(\sigma)}=\phi(s_{\overline w_2,\overline w_1})\in \phi(\cl A_D)$.
This holds for each $\sigma\in \overline{\fk W}_m$.
Therefore
$$
g= \displaystyle \sum_{\sigma \in \overline{\fk W}_m} g\chi_{\Omega(\sigma)}\in \phi(\cl A_D).
$$
This shows that the range of $\phi$ contains $\Gamma$, as required.
\end{proof}

We now prove that conditions (H2) and (H3) are satisfied.
To verify condition (H2), it is enough to show that if $a, b \in A$, then there is a directed segment
$\sigma$ such that
\begin{equation}\label{irred}
a=\Gamma o(\sigma),\quad b= \Gamma t(\sigma).
\end{equation}
Let $a=\Gamma \sigma_1$, $b=\Gamma \tau_2$, where $\sigma_1, \tau_2\in\fk S_{k+1}$. By \cite[Proposition 1 (iii)]{ch}, there is a $\Gamma$-periodic geodesic $\gamma$ containing $\tau_2$. By definition, this means that there is a subgroup of $\Gamma$ which leaves the geodesic $\gamma$ invariant and acts upon it by translation. Choose $\omega\in \partial\Delta$ to be the boundary point of $\gamma$ with
$$\tau_2\subset[o(\tau_2), \omega).$$
Since the action of $\Gamma$ on $\partial\Delta$ is minimal, there exists $g\in \Gamma$ such that $g\omega\in \Omega(\sigma_1)$.  The geodesic $g\gamma$ is $\Gamma$-periodic. Therefore the  semi-geodesic $[go(\tau_2), g\omega)$ contains infinitely many directed segments $\sigma_2$  which are $\Gamma$-translates of $\tau_2$. Choose such a segment $\sigma_2$ far enough away from $g\tau_2$ so that $\sigma_2\in \Omega(\sigma_1)$. Let $\sigma$ be the directed segment with $o(\sigma)=\sigma_1$ and  $t(\sigma)=\sigma_2$. Then (\ref{irred}) is satisfied.

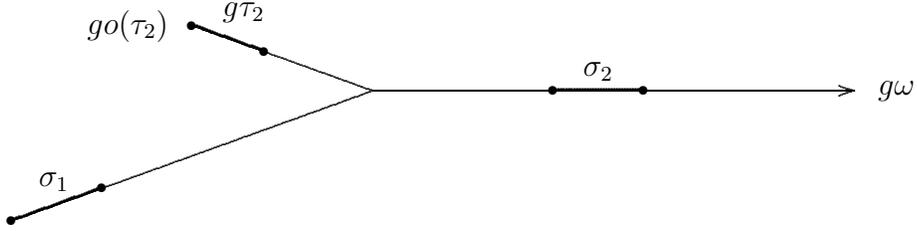
\centerline{
\beginpicture
\setcoordinatesystem units <0.8cm,0.866cm>   
\setplotarea x from -3 to 10, y from -2  to 2         
\putrule from 2 0 to 10 0
\setlinear
\plot  -1 1  2 0 /
\plot  -4 -2  2 0 /
\put{${g\omega}$}[l]  at  10.4  0
\put{$_\bullet$}  at  -1 1
\put{$_\bullet$}  at  0.2 0.6
\put{$g\tau_2$}[r]  at  0.2 1.2
\put{$go(\tau_2)$}[r]  at  -1.4 1
\put{$_\bullet$}  at  -4 -2
\put{$_\bullet$}  at  -2.5 -1.5
\put{$\sigma_1$}  at  -3.3 -1.4
\put{$_\bullet$}  at  5 0
\put{$_\bullet$}  at  6.5 0
\put{$\sigma_2$}  at  5.75 0.3
\arrow <6pt> [.3,.67] from  9.8 0 to  10 0
\setplotsymbol({$\cdot$}) \plotsymbolspacing=1pt
\plot  -1 1   0.2 0.6 /
\plot  5 0   6.5 0 /
\plot  -4 -2  -2.5 -1.5 /
\endpicture
}

\bigskip

To prove that condition (H3) holds, let $p>0$. Since $\Delta$ has more than two ends, there exist vertices of $\Delta$ which have degree greater than $2$. Let $\sigma=(s_0, s_1, \dots , s_{p+k})\in \fk S_{p+k}$ be a directed segment whose final vertex $s_{p+k}$ has degree greater than $2$.
Extend $\sigma$ to two different segments:
$$(s_0, s_1, \dots , s_{p+k}, t),\, (s_0, s_1, \dots , s_{p+k}, t')\in \fk S_{p+k+1}.$$

\centerline{
\beginpicture
\setcoordinatesystem units <0.5cm,0.866cm>   
\setplotarea x from -10 to 2, y from -2  to 2         
\putrule from 0 0 to -10 0
\setlinear
\plot 2 1  0 0  2 -1 /
\put{$_\bullet$}  at  0 0
\put{$_\bullet$}  at  -10 0
\put{$_\bullet$}  at  2 1
\put{$_\bullet$}  at  2 -1
\put{$s_0$}  at  -11 0
\put{$s_{p+k+1}$}  at  1.35  0
\put{$t$}  at  3 1
\put{$t'$}  at  3 -1
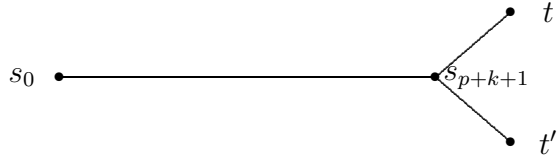
\endpicture
}

Let $w=\alpha(\Gamma (s_0, \dots , s_{p+k}, t))$ and $w'=\alpha(\Gamma (s_0, \dots , s_{p+k}, t'))$. Then $w, w'\in W_p$, $o(w)=o(w')$ and
$t(w)\ne t(w')$ since $\Gamma$ acts freely on $\fk S_k$ and $\alpha$ is injective. Therefore at least one of the words $w, w'$ is not $p$-periodic.

\section{K-theory and examples}

The Standing Hypotheses (1)--(4) remain in force. Thus the algebra $\cl A_\Gamma$ is isomorphic to the Cuntz-Krieger algebra $\cl A_D$ associated with the alphabet $A$, the set of decorations $D$ and transition matrix $M$.

The simple Cuntz-Krieger algebras $\cl A_D$ are purely infinite, nuclear and satisfy the Universal Coefficient Theorem \cite[Remark 6.5]{rs1}. They are therefore classified by their K-theory \cite{k}.
It is convenient to consider the related algebra $\cl A_A$ which is stably isomorphic to $\cl A_D$.
Recall from Remark \ref{OA} that $\cl A_A$ is isomorphic to the algebra $\cl O_{M^t}$.
The groups $K_0(\cl A_A)$, $K_1(\cl A_A)$ are the Bowen-Franks invariants of flow equivalence for a certain subshift associated with $(\Gamma, \Delta)$. More precisely, according to \cite[Proposition 3.1]{c} the  group $K_0(\cl A_A)$ is isomorphic to the abelian group
\begin{equation}\label{K}
\cl G_\Gamma=\left\langle A\, \left|\, a=\sum_{b\in A} M(a,b)b,\, a\in A\right.\right\rangle\,.
\end{equation}
Note that, as the notation suggests, $\cl G_\Gamma$ depends only on $\Gamma$, by Remark \ref{qi}.
Also $K_1(\cl A_A)$ is the torsion free part of $\cl G_\Gamma$. Therefore $\cl A_A$ is classified up to stable isomorphism by the group $\cl G_\Gamma$. Since the algebra $\cl A_D$ is stably isomorphic to $\cl A_A$ \cite[Corollary 5.15]{rs1}, we obtain the following result.

\begin{theorem}\label{Ktheorem}
Under the Standing Hypotheses (1)--(4),
$K_0(\cl A_\Gamma)\cong \cl G_\Gamma.$
\end{theorem}

To completely classify $\cl A_A$ up to isomorphism, we need to identify the class $[\id]$ of the identity idempotent in $K_0(\cl A_A)$ \cite{k}. By Remark \ref{id}, this class corresponds to the element
\begin{equation}\label{identity}
  \varepsilon =\displaystyle\sum_{a\in A}a \in \cl G_\Gamma.
\end{equation}
Here are explicit calculations in the case where $\Gamma$ is a free product of finite cyclic groups.

\subsection{Example:}\label{ex1} The group $\Gamma=\bb Z_{l+1}*\bb Z_{m+1}$ acts on its Bass-Serre tree \cite[I.4]{ser} with an edge $y=[P,Q]$ as its fundamental domain. The stabilizer $\Gamma_P$ of $P$ is isomorphic to $\bb Z_{l+1}$ and the stabilizer $\Gamma_Q$ of $Q$ is isomorphic to $\bb Z_{m+1}$.

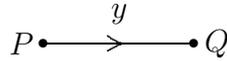
\centerline{
\beginpicture
\setcoordinatesystem units <.5cm,1cm>   
\setplotarea x from -4 to 4, y from -1  to 1         
\putrule from -2 0 to 2 0
\arrow <6pt> [.4, 1] from  0 0  to 0.1 0
\put{$_\bullet$}  at  2 0
\put{$_\bullet$}  at  -2 0
\put{$P$}  at  -2.6 0
\put{$Q$}  at  2.6 0
\put{$y$}  at  0 0.4
\endpicture
}

By construction, $\Gamma$ acts freely and transitively on the geometric edges of $\Delta$. In other words, $\Gamma$ is 1-acylindrical.
The theory applies, with $k=1$, and the alphabet $A$ is the set of $\Gamma$-orbits of directed segments of length $2$ in $\Delta$.

Let $A_1=\Gamma_P-\{1\}$ and $A_2=\Gamma_Q-\{1\}$, so that $|A_1|=l$, $|A_2|=m$.  Each directed segment of length $2$ in $\Delta$ lies in the $\Gamma$-orbit of one of the directed segments $[P, a_2P]$, $[Q, a_1Q]$, for some  $a_1\in A_1$, $a_2\in A_2$.
Let $\hat a_2=\Gamma[P, a_2P]$, $\hat a_1=\Gamma[Q, a_1Q]$ be the corresponding elements of $A_2$, $A_1$.

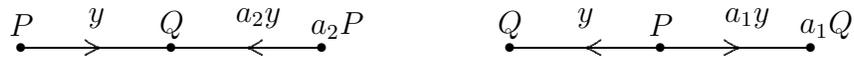
\centerline{
\beginpicture
\setcoordinatesystem units <.5cm,1cm>   
\setplotarea x from -2 to 6, y from -1  to 1         
\putrule from -2 0 to 6 0
\arrow <6pt> [.4, 1] from  0 0  to 0.1 0
\arrow <6pt> [.4, 1] from  4.1 0  to 4 0
\put{$_\bullet$}  at  6 0
\put{$_\bullet$}  at  2 0
\put{$_\bullet$}  at  -2 0
\put{$P$}  at  -2 .3
\put{$Q$}  at  2 .3
\put{$a_2P$}  at  6.4 .3
\put{$y$}  at  0 0.4
\put{$a_2y$}  at  4.3  .4
\setcoordinatesystem units <.5cm,1cm>  point at -13 0   
\setplotarea x from -2 to 6, y from -1  to 1        
\putrule from -2 0 to 6 0
\arrow <6pt> [.4, 1] from  0.1 0  to 0 0
\arrow <6pt> [.4, 1] from  4 0  to 4.1 0
\put{$_\bullet$}  at  6 0
\put{$_\bullet$}  at  2 0
\put{$_\bullet$}  at  -2 0
\put{$Q$}  at  -2 .3
\put{$P$}  at  2 .3
\put{$a_1Q$}  at  6.4 .3
\put{$y$}  at  0 0.4
\put{$a_1y$}  at  4.3  .4
\endpicture
}
The map $a\mapsto \hat a$ is a bijection from $A_1\cup A_2$ onto $A$.
The $\{0,1\}$-matrix $M$ is defined by
$M(\hat a,\hat b)=1 \Longleftrightarrow$ either $a\in A_1, b\in A_2$ or $b\in A_1, a\in A_2$.

\centerline{
\beginpicture
\setcoordinatesystem units <.5cm,1cm>   
\setplotarea x from -2 to 6, y from -1.5  to 1         
\putrule from -2 0 to 10 0
\put{$_\bullet$}  at  10 0
\put{$_\bullet$}  at  6 0
\put{$_\bullet$}  at  2 0
\put{$_\bullet$}  at  -2 0
\put{$P$}  at  -2 .3
\put{$Q=a_2Q$}  at  2 .3
\put{$a_2a_1Q$}  at  10.4 .3
\put{$a_2P$}  at  6.4 .3
\put{The condition $M(\hat a_1, \hat a_2)=1$.} at 4 -0.8
\endpicture
}

By Theorem \ref{Ktheorem},
\begin{equation}\label{Ka}
K_0(\cl A_\Gamma)=\left\langle \hat A_1\cup\hat A_2 \, \left|\, \hat a=\sum_{b\in A_{3-j}}\hat b,\quad a\in A_j, \quad j=1,2\right.\right\rangle\,.
\end{equation}
The relations on the right side of (\ref{Ka}) show that all the generators of $\hat A_1$ are equal and all the generators of $\hat A_2$ are equal. Therefore
\begin{equation*}\label{Kab}
K_0(\cl A_\Gamma)=\langle \hat a_1, \hat a_2\, \left|\right. \hat a_2=l\hat a_1, \hat a_1=m \hat a_2\rangle
=\langle \hat a_1\, \left|\right. \hat a_1=lm\hat a_1\rangle
=\bb Z_{lm-1}.
\end{equation*}
Recall that the classical Cuntz algebra $\cl O_n$ is generated by $n$ isometries whose range projections sum to the identity operator \cite{c0}. Now $K_0(\cl O_n)= \bb Z_{n-1}$ \cite{ca}. It follows from the classification theorem \cite{k} that $\cl A_\Gamma$ is stably isomorphic to $\cl O_{lm}$.

In order to classify $\cl A_\Gamma$ up to isomorphism, the class $[\id]$ of the identity idempotent in $K_0(\cl A_\Gamma)$ must be identified. Now $[\id]$ corresponds to the element $(l+1)\hat a_1 =\hat a_1+\hat a_2=(m+1)\hat a_2$.
However it is known \cite{ca} that
$$(K_0(M_k\otimes \cl O_n), [\id]) \cong (\bb Z_{n-1}, k).$$
This proves that $\cl A_\Gamma \cong M_{l+1}\otimes \cl O_{lm}\cong M_{m+1}\otimes \cl O_{lm}$.

\subsection{Example:} More generally, a free product of finite groups $\Gamma=\Gamma_1*\Gamma_2*\dots*\Gamma_n$ acts on its Bass-Serre tree. The fundamental domain is a tree consisting of $n$ edges $y_i$ with terminal vertex $Q_i$ emanating from a common vertex $P$. The stabilizer of $P$ and of each of the edges $y_j$ is trivial, and the stabilizer of $Q_j$ is isomorphic to $\Gamma_j$.

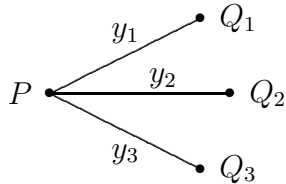
\centerline{
\beginpicture
\setcoordinatesystem units <1cm,1cm>   
\setplotarea x from -2 to 2, y from -2  to 2         
\putrule from 0 0 to 2.4 0
\setlinear
\plot 2 1  0 0  2 -1 /
\put{$_\bullet$}  at  0 0
\put{$_\bullet$}  at  2 1
\put{$_\bullet$}  at  2 -1
\put{$_\bullet$}  at  2.4 0
\put{$Q_1$}  at  2.5 1
\put{$Q_2$}  at  2.9 0
\put{$Q_3$}  at  2.5 -1
\put{$P$} at -0.4 0
\put{$y_2$} at 1.5 0.2
\put{$y_1$} at 1 0.8
\put{$y_3$} at 1 -0.8
\endpicture
}

The group $\Gamma$ is 1-acylindrical and acts freely (but not transitively, in contrast to Example \ref{ex1}) on the set of edges of $\Delta$, with finitely many orbits.

Let $A_i=\Gamma_i-\{1\}$ and $\gamma_i=|A_i|$, $1\le i\le n$.  Each directed segment of length $2$ in $\Delta$ lies in the $\Gamma$-orbit
of a directed segment of the form $(P, Q_i, aP)$, $a\in A_i$, or $(Q_j, P, Q_k)$, $j\ne k$, as illustrated below.

\centerline{
\beginpicture
\setcoordinatesystem units <.5cm,1cm>   
\setplotarea x from -2 to 6, y from -1  to 1         
\putrule from -2 0 to 6 0
\put{$_\bullet$}  at  6 0
\put{$_\bullet$}  at  2 0
\put{$_\bullet$}  at  -2 0
\put{$P$}  at  -2 .3
\put{$Q_i$}  at  2 .3
\put{$y_i$}  at  0 0.2
\put{$ay_i$}  at  4.3  .2
\put{$aP$}  at  6.2 .3
\setcoordinatesystem units <.5cm,1cm>  point at -13 0 
\setplotarea x from -2 to 6, y from -1  to 1         
\putrule from -2 0 to 6 0
\put{$_\bullet$}  at  6 0
\put{$_\bullet$}  at  2 0
\put{$_\bullet$}  at  -2 0
\put{$Q_j$}  at  -2 .3
\put{$P$}  at  2 .3
\put{$y_j$}  at  0 0.2
\put{$y_k$}  at  4.3  .2
\put{$Q_k$}  at  6 .3
\endpicture
}

Let $\hat A_i = \{\hat a=\Gamma(P, Q_i, aP) :  a\in A_i\}$,
and let $\hat B = \{\hat b_{jk}= \Gamma(Q_j, P, Q_k) :  j\ne k\}$. Then by Theorem \ref{Ktheorem}
\begin{equation}\label{Ka2}
K_0(\cl A_\Gamma)=\left\langle \bigcup_i\hat A_i \cup\hat B \, \left|\, \hat a=\sum_{j\ne i} \hat b_{ij} \quad (\hat a\in\hat A_i),\quad
\hat b_{jk}=\sum_{\hat a\in \hat A_k} \hat a
 \right.\right\rangle.
\end{equation}

 The relations on the right side of (\ref{Ka2}) show that $\hat b_{jk}$ depends only on $k$. Therefore, for each $i$, all the generators  in $\hat A_i$ are equal. It follows that
\begin{equation}\label{Ka3}
K_0(\cl A_\Gamma)=\langle \hat a_i\, \left|\right. \hat a_i = \sum_{j\ne i} \gamma_j\hat a_j\rangle.
\end{equation}
It is easy to see from (\ref{Ka3}) that $K_0(\cl A_\Gamma)$ is a torsion group and therefore that $K_1(\cl A_\Gamma)=0$. In other words, the unitary group of $\cl A_\Gamma$ is connected \cite{ca}.
If all the groups $\Gamma_i$ have the same order, say $\gamma+1$, and $\delta=\gamma(n-1)-1$, then (\ref{Ka3}) simplifies to
\begin{equation*}
K_0(\cl A_\Gamma)=\bb Z_{(\gamma+1)\delta}\oplus(\bb Z_{\gamma+1})^{n-2}
\end{equation*}
with canonical generators $\hat a_1, \hat a_2- \hat a_1, \hat a_3- \hat a_1, \dots, \hat a_{n-1}- \hat a_1$.

\end{document}